\documentclass[12pt]{amsart}

\usepackage{amsmath,amssymb,amscd,amsfonts,amsthm,verbatim}
\usepackage[all,cmtip]{xy}
\usepackage{paralist}
\usepackage[mathscr]{eucal}
\usepackage{color}
\usepackage{pdfsync}
\usepackage[]{fontenc}
\usepackage{enumerate}
\usepackage{pgfplots}

\usepackage[pdftex]{hyperref}
\hypersetup{citecolor=blue,linktocpage}

\usepackage[colorinlistoftodos]{todonotes}

\newtheorem{thm}{Theorem}[section]
\newtheorem{lem}[thm]{Lemma}
\newtheorem{prop}[thm]{Proposition}

\newtheorem{cor}[thm]{Corollary}

\newtheorem{assu-nota}[thm]{Assumption--Notation}

\newtheorem*{thm*}{Theorem}

\theoremstyle{definition}
\newtheorem{defn}[thm]{Definition}
\newtheorem{rem}[thm]{Remark}
\newtheorem{ex}[thm]{Example}

\newcommand{\R}{\mathbb R}
\newcommand{\F}{\mathbb F}

\DeclareMathOperator{\alb}{alb}

\DeclareMathOperator{\vol}{vol}

\numberwithin{equation}{section}

\title{Slope Inequalities for higher dimensional irregular fibrations}
\author{Miguel \'Angel Barja}
\address{Miguel \'Angel Barja\\Departament de Matem\`atiques\\Universitat Polit\`ecnica de Catalunya\\ Institut de Matemàtiques de la UPC-BarcelonaTech (IMTech) \\Avda. Diagonal 647\\08028 Barcelona\\Spain}
\email{miguel.angel.barja@upc.edu}

\thanks{The author was supported by MINECO PID2019-103849GB-I00 ``Geometr{\'\i}a Álgebra, Topolog{\'\i}a y aplicaciones multidisciplinares" and by Generalitat de Catalunya SGR2018.}

\begin{document}
\begin{abstract}
We prove the equivalence between Clifford-Severi inequalities for good classes of varieties of maximal Albanese dimension and Slope Inequalities for such varieties fibred over curves. This provides a big set of new Slope Inequalities and characterizes the limit cases. It also gives a machinery to automatically obtain other higher dimensional Slope and Clifford-Severi inequalities from inequalities in low dimension. For this, we construct a continuous version of Xiao's method for irregular fibrations.

\end{abstract}

\maketitle
\setcounter{tocdepth}{1}
\tableofcontents
\section{Introduction}

\bigskip

Given a $n$-dimensional variety $X$ and a big line bundle $L$ it is a natural question to ask for lower bounds of the volume of $L$ in terms of the global sections $h^0(X,L)$. These are the so called {\it Noether-type inequalities}. Clifford's lemma for curves or Noether's inequality for surfaces are examples of this. If the variety is fibred over a curve $f:X\longrightarrow B$, another natural question is to find lower bounds for the {\it slope} of $L$ with respect to $f$:

$$
s(f,L)=\frac{{\rm vol}(L)}{\deg f_*L}
$$

\noindent when $\deg f_*L>0$. We call them {\it Slope Inequalities}. The slope inequality of Xiao, Cornalba-Harris and Moriwaki is the standard example.

Both kind of inequalities are closely related but in general there is not a clear implication between them. When the variety $X$ is irregular, the relation is clearer. For example, the celebrated proof of Pardini of the Severi conjecture, deduces the Severi inequality from Xiao's slope inequality. When we consider the extra data of a nontrivial map $a: X \longrightarrow A$ to an abelian variety, two modifications of the former invariants come into the scene: the {\it continuous rank} $h^0_a(X,L)$ (see \cite{B} and \cite{BPS3}) and the {\it positive continuous degree} $\deg^+_af_*L$ that we introduce in Section 2.

In the context of irregular varieties, for big line bundles, we can consider

$$
\lambda (L,a)=\frac{{\rm vol} (L)}{h^0_a(X,L)}
$$

\noindent and {\it Clifford-Severi} inequalities are lower bounds for this ratio.

The main example in dimension 2 is the original work of Pardini on the Severi conjecture (\cite{P}), which gives the name to the bounds. Also for surfaces there are further results by Mendes-Lopes-Pardini (\cite{MLP}) and by Lu and Zuo ([\cite{LZ}). In the last ten years, several higher dimensional Clifford-Severi bounds have been obtained. The first results were found independently by Zhang in \cite{Z} (for $L=K_X$ and maximal Albanese dimension) and the author in \cite{B} (for general $L$ and $a$). Finally, in a series of works (\cite{BPS2}, \cite{BPS3} and \cite{BPS4}; see also \cite{S}), Pardini, Stoppino and the author introduce new techniques to study linear series on irregular varieties that, as a byproduct, produce a huge series on higher dimensional Clifford-Severi inequalities for varieties of maximal Albenese dimension. Moreover, in the two main cases, pairs $(X,L)$ verifying the limit equality can be classified. Clifford Severi inequalities for varieties of non-maximal Albanese dimension where also obtained by the author in \cite{B} and later on were extended by Jiang (\cite{J}) considering the geometry of the fibre of the Albanese map. In the case of Albanese dimension $n-1$, Zhang (\cite{Z2}) obtained stronger inequalities and also characterized the limit cases.

In \cite{B}, the author proposes as a conjecture whether similar {\it slope inequalities} hold for irregularly fibred varieties, and prove them in the case of low base genus. Recently, Hu and Zhang (\cite{HZ}) obtain these higher dimensional {\it slope inequalities} for fibrations of maximal Albanese dimension and $L=K_f$. The arguments used there are variants of those in \cite{BPS3} and \cite{Z} by a direct computation using classical Xiao's method. They also give properties of fibrations that attain the equality in the main case.

\smallskip

The aim of this paper is threefold. First of all, we will show that {\it Clifford-Severi} inequalities and {\it Slope inequalities} for good classes of varieties of maximal Albanese dimension are equivalent. More concretely, fix a {\it good} class ${\mathcal F}$ of triplets $(X,L,a)$ (see Section 4 for a definition and examples). As an instance, one important case is $X$ to be of maximal $a$-dimension and $L$ to be numerically subcanonical.

Denote by ${\lambda}_{\mathcal F}(n)$ to be the infimum of Clifford-Severi slopes for triplets of dimension $n$ in the class ${\mathcal F}$. Similarly, we denote by $s_{\mathcal F}(n)$ to be the infimum of the slopes $s(f,L,a)=\frac{{\rm vol}(L)}{\deg_a^+ f_*L}$
 for irregular fibrations  $f: X\longrightarrow B$, $n=\dim X$, with general fibre $F$ such that $(F,L_{|F},a_{|F})$ is in ${\mathcal F}$.
For example, we can see the Clifford-Severi inequalities given in Section 5 of \cite{BPS3} as a set of inequalities proved by induction on the dimension in a way that can be rephrased as:

$$
{\lambda}_{\mathcal F}(n)\geq \,n\,{\lambda}_{\mathcal F}(n-1).
$$

In this paper we prove that the above inequality can be split (see Theorem \ref{metateorema}), given a different and extremely simple proof:

\bigskip

\noindent {\bf Main Theorem}. \label{main} ${\lambda}_{\mathcal F}(n)\geq \, s_{\mathcal F}(n)\geq \,n\, {\lambda}_{\mathcal F}(n-1).$


\bigskip

A first consequence of this Theorem is that using the whole set of the known Clifford-Severi inequalities (e.g. \cite{B}, \cite{BPS3}) we can deduce similar slopes inequalities. But, even more importantly, this theorem constitutes a machinery to automatically obtain higher dimensional Slope and Clifford-Severi inequalities from low dimensional results, via a zig-zag argument. In particular, it provides an alternative and simple way to obtain all the aforementioned known Clifford-Severi inequalities (and the corresponding Slope inequalities) just proving them in dimension 1.

\bigskip

\noindent{\bf Corollary A}. {\it If, for a good class of varieties of maximal Albanese dimension ${\mathcal F}$, we have $\lambda_{\mathcal F}(m_0)\geq \lambda_{m_0}$ for a certain low dimension $m_0$, then, for triplets in ${\mathcal F}$ we have

$$
\lambda_{\mathcal F}(n) \geq \frac{n!}{m_0!}\lambda_{m_0}
$$

\noindent and for fibred varieties of dimension $n$ with fibres in the class ${\mathcal F}$ we have

$$
s_{\mathcal F}(n)\geq \frac{n!}{{m_0}!}\lambda_{m_0}.
$$

}

\bigskip

Our second objective is to study the equality in the Slope inequalities. Consider a fibration $f:X \longrightarrow B$, $n=\dim X$,  with general fibre $F$ such that $(F,L_{|F},a_{|F})$ belongs to a good class ${\mathcal F}$. Assume that the slope of this fibration attains the minimal slope $n\,{\lambda}_{\mathcal F}(n-1)$. Then we prove that the fibres verify a Clifford-Severi {\it equality} and many other properties hold (see theorem \ref{metateorema} and Corolaries \ref{cor1} and \ref{cor2}). These results generalize and extend those of \cite{HZ} and can also be viewed as a huge extension of a theorem of Xiao (\cite{X}, Theorem 3) on irregularly fibred surfaces of slope 4. Here we have a typical example of slope inequality and equality (see Corollary \ref{corolario}):

\bigskip

\noindent{\bf Corollary B}.
{\it Let $f:X \longrightarrow B$ a relatively minimal fibration of maximal Albanese dimension. Let $F$ be a general fibre. Then:
\begin{itemize}
\item[(i)] $K_f^n\geq 2\,n! \chi_f.$
\item[(ii)] If equality holds, then
    \begin{itemize}
      \item [$(1)$] $K_F^{n-1}=2\,(n-1)! \chi(\omega_F).$
      \item [$(2)$] For general $\alpha \in {\rm Pic}^0X$, the vector bundle $f_*(\omega_f\otimes \alpha)$ is nef and semistable.
      \item [$(3)$] If $\chi(\omega_F)>0$ and $n\geq 3$ then $q(F)=n-1$, $q(X)=b+q(F)$ and the map $(alb_X)_{|F}:F\longrightarrow K$ is, birationally, a double cover composed with an étale map on $K$.
    \end{itemize}

\end{itemize}
}

\bigskip

Let us come back to the statement of the Main Theorem. The first inequality states that a Slope Inequality for a class ${\mathcal F}$ implies the same Clifford-Severi inequality for this class and the proof follows by applying an adapted version of Pardini's trick.

The second inequality in Main Theorem follows easily from a continuous version of Xiao's method. The construction of this {\it continuous Xiao's method} is the third objective of this paper. We must emphasize here that using this continuous version, the argument to prove this inequality is extremely simple, compared with that of \cite{HZ} and does not depend on the concrete class ${\mathcal F}$.

\bigskip

The paper is organized as follows. In Section 2 we introduce the (positive) continuous degree of a line bundle on an irregular variety with respect to an irregular fibration and a map to an abelian variety, and study its properties. This is the right invariant to study Slope Inequalities in the context of irregular fibration. There we also compare this invariant with the continuous rank as this information turns out to be fundamental to compare Clifford-Severi and Slope inequalities.

In Section 3 we construct the whole theory of the {\it continuous Xiao's method}.

In Section 4 we state and prove the Main Theorem for good classes of maximal $a$-dimension varieties, we deduce a set of new slope inequalities and the characterization of the limit cases.

\bigskip
\bigskip

\noindent  {\underline{Notations and Conventions}} We work over $\mathbb{C}$. Varieties are projective and smooth unless otherwise stated. We will use notation of divisor or line bundles indistinctly.

\bigskip

\noindent {\underline{Acknowledgements}} The author thanks Lidia Stoppino and Rita Pardini for extremely useful discussions on this topic along the last years. The main ideas of the paper were developed during a short stay at the Department of Mathematics of Pavia University. The author thanks Gian Pietro Pirola, Lidia Stoppino and all the members of the Algebraic Geometry group of this department for their warm hospitality.

\bigskip
\bigskip

\section{The continuous degree}

In this section we will define and study the {\it (positive) continuous degree} of a line bundle on an irregularly fibred variety.

\begin{assu-nota}\label{assumptions}

Let $X$ be a complex, projective variety of dimension $n$ and $a:X\longrightarrow A$ a nontrivial map to an abelian variety. We will assume that $a$ is {\it strongly generating}, i.e., that $a(X)$ generates $A$ and the pullback map $a^*:{\widehat A}\rightarrow {\rm Pic}^0(X)$ is injective.

Consider a fibration $f: X\longrightarrow B$ onto an smooth curve of genus $b$ with general smooth fibre $F$. According to \cite{Beau}, the natural restriction map $i^*: {\rm Pic}^0(X) \longrightarrow {\rm Pic}^0 (F)$ has kernel which is a finite extension of $f^*({\rm Pic}^0 (B))$. Consider this situation restricted to $a^*({\widehat A})$, and let ${\widehat K}:=i^*(a^* ({\widehat A}))={a_F}^*({\widehat A})$ in ${\rm Pic}^0(F)$. Consider the subtorus $${\widehat T}=(a^*)^{-1}(f^*{\rm Pic}^0(B))\subseteq {\widehat A}$$

\noindent and the corresponding projection $A\longrightarrow T$. Altogether this fits into a commutative diagram

\begin{equation}\label{basicdiagram}
\xymatrix {F\ar[d]^{i}\ar[r]^{a_{|F}}& K \ar[d]\\X \ar[d]^{f}\ar[r]^{a} & A \ar[d] \\B\ar[r]^{a_B} &T}
\end{equation}
\bigskip

Observe that $a(F)$ generates $K$, but the map $a_{|F}:F\longrightarrow K$ needs not to be strongly generating since the map ${\widehat A}\rightarrow {\widehat K}$ may not have connected kernel.

We denote by $q=\dim A\geq 1$, ${\overline b}=\dim T\geq 0$ and $q_f=\dim K=q-{\overline b}$.

We will always assume that the {\it fibration $f$ is irregular}, i.e.,${\rm dim} K=q_f\geq 1$.

We say that $f$ is of $a$-dimension $k$ if $k=\dim a(F)$. We say that $f$ is of maximal $a$-dimension if $k=n-1$. Observe that if $X$ if of maximal $a$-dimension, then so is $f$ but the viceversa does not hold in general if ${\overline b}=0$ (cf. \cite{HZ}). When $a={\rm alb}_X$ these correspond to the usual concepts of maximal Albanese dimensions and by the universal property of the Albanese variety we have that $b={\overline b}$.

\end{assu-nota}

\bigskip
\bigskip

Consider the family $\{f_*(L\otimes \alpha)\}_{{\alpha}\in {\widehat A}}$ of vector bundles on $B$. It is flat over a non-empty open set $U\subseteq {\widehat A}$. Restricting to an even smaller open set $U_0$, we have that the length, ranks and degrees of their Harder-Narashiman ($HN$) filtrations are constant. Moreover, they are upper semicontinuous on $U$(\cite{Sh}). For $\alpha \in U_0$ denote

$$
0={\mathcal E}^{\alpha}_0\subset {\mathcal E}^{\alpha}_1\subset ...\subset  {\mathcal E}^{\alpha}_l={\mathcal E}^{\alpha}=f_*(L\otimes \alpha)
$$

\bigskip

\noindent its $HN$ filtration and the slopes $\mu_1>\mu_2>...>\mu_l$, ($\mu_i=\mu({\mathcal E}^{\alpha} _i/{\mathcal E}^{\alpha} _{i-1}$), degrees $d_i=\deg {\mathcal E}^{\alpha}_i$ and ranks $r_i={\rm rk}{\mathcal E}^{\alpha}_i$.

\begin{defn} For $\alpha \in U_0$, let $m$ ($0\leq m\leq l$) be the biggest index such that $\mu_m\geq 0$ and denote the {\it positive part} as $f_*(L\otimes \alpha)^+={\mathcal E}^{\alpha}_m $ (by abuse of notation, se say $m=0$ if $\mu_i<0$ for all $i$).
\end{defn}

\begin{lem}\label{remarkpositivepart}

\begin{itemize}
\item [(i)] $f_*(L\otimes \alpha)^+$ is the biggest nef subbundle of $f_*(L\otimes \alpha)$, i.e., if ${\mathcal G}\subseteq f_*(L\otimes \alpha)$ is nef, then  ${\mathcal G}\subseteq f_*(L\otimes \alpha)^+$.
\item [(ii)] For $\alpha \in U_0$, $h^0_a(X,L)=h^0(B,{\mathcal E}^{\alpha}_m)$. In particular, if $h^0_a(X,L)\neq 0$, then $m\geq 1$.

\end{itemize}
\end{lem}

\begin{proof} (i) Assume $m<l$, otherwise the result is trivial. Consider the projection $\pi: {\mathcal E}^{\alpha} \rightarrow {\mathcal E}^{\alpha}/{\mathcal E}^{\alpha}_{m}$. The vector bundle $\pi ({\mathcal G})$ is either 0 or it has non-negative slope, since ${\mathcal G}$ is nef. It is a subsheaf of ${\mathcal E}^{\alpha}/{\mathcal E}^{\alpha}_{m}$ which has $\mu_+({\mathcal E}^{\alpha}/{\mathcal E}^{\alpha}_{m})=\mu_{m+1}<0$. Hence $\pi ({\mathcal G})=0$, and so ${\mathcal G}\subseteq {\mathcal E}_m$.

\smallskip
\noindent (ii) We have that $h^0_a(X,L)=h^0(X,L\otimes \alpha)=h^0(B,f_*(L\otimes \alpha))$. Now, observe that for $i>m$ we have $\mu_i<0$ and so $h^0(B,{\mathcal E}_i/{\mathcal E}_{i-1})=0$. By induction on $i>m$ we obtain $h^0(B,{\mathcal E}_m^{\alpha})=h^0(B,{\mathcal E}_i^{\alpha})$.
\end{proof}

\begin{defn} The {\it continuous degree} ${\deg}_af_*L$ and {\it continuous rank} ${\rm rk}_a f_*L$ are those of $f_*(L\otimes \alpha)$ for $\alpha \in U_0$ as above.

Analogously, we define the {\it continuous positive degree} ${\deg}^+_af_*L$ and {\it continuous positive rank} ${\rm rk}^+_a f_*L$ to be those of $f_*(L\otimes \alpha)^+$ for $\alpha \in U_0$.

\end{defn}

\medskip

\begin{ex} When $L=K_f$ the vector bundle $f_*K_f$ is always nef. If $f$ is of maximal $a-$dimension, we can apply Generic Vanishing on $F$ and obtain that $\deg^+ _a f_*(K_f)=\deg _a f_*(K_f)=\chi _f$ (see \cite{HP}, \cite{HZ}).
\end{ex}
\medskip

\begin{rem}{\label{torsion}} Let $D={\overline A}\setminus U_0$. Fix $\alpha _0 \notin \bigcup_{\tau \in {\widehat A}_{\rm tors}}(\tau+D)$ (so it is very general). Then for any $\tau \in {\widehat A}_{\rm tors}$ we have that $\alpha_0+\tau\in U_0$ and so all the above continuous invariants correspond to that of $f_*(L_0\otimes \tau)$, where $L_0=L\otimes \alpha_0$.
\end{rem}
\medskip

\begin{rem}\label{elecciond}
We can consider the basic diagram (\ref{basicdiagram}) after a base change via a multiplication map (see \cite{BPS3}). For this, fix an integer $d$ and consider the corresponding multiplication maps by $d$ on $A$, $K$ and $T$ respectively. Consider the corresponding base changes and we obtain

$$\xymatrix {F^{(d)}\ar[d]^{i^{(d)}}\ar[r]^{{(a_d)}_{|F^{(d)}}}& K \ar[d]\\X^{(d)} \ar[d]^{f^{(d)}}\ar[r]^{a_d} & A \ar[d] \\B^{(d)}\ar[r]^{a_B^{(d)}} &T}$$

\bigskip

We will use the notation $\nu_d$ for the multiplication map at the level $B,X,F,T,A$ and $K$ if no confusion arises. Observe that these maps are ètale of degree $d^{2\overline b}$ on $T$ and $B$, of degree $d^{2q}$ on $A$ and $X$ and of degree $d^{2q_f}$ on $K$ and $F$.

If $d$ is coprime with the cardinality of the finite part of ${\rm Ker} a_{F}^*:{\widehat A}\rightarrow {\widehat K}$, we have we have an exact sequence of the subgroups of $d$-torsion $$0\rightarrow {\widehat T}_d\rightarrow {\widehat A}_d\rightarrow {\widehat K}_d\rightarrow 0.$$

In particular, $F^{(d)}$ is connected and so $f^{(d)}$ is again a fibration. In what follows, we will consider only multiplication maps by such $d$.

\end{rem}

\medskip

\medskip

\begin{lem} Consider the notations and hypotheses of (\ref{assumptions}) and (\ref{elecciond}). Let $L^{(d)}:=\nu _d^* L$. Then we have

\begin{itemize}
\item [(i)] $\deg _{a_d}f^{(d)}_*L^{(d)}=d^{2q}{\deg _a f_*L}$ and $\deg^+ _{a_d}f^{(d)}_*L^{(d)}=d^{2q}{\deg^+ _a f_*L}$

\medskip

\item [(ii)] ${\rm rk} _{a_d}f^{(d)}_*L^{(d)}=d^{2q_f}{\rm rk} _a f_*L$ and ${\rm rk}^+ _{a_d}f^{(d)}_*L^{(d)}=d^{2q_f}{\rm rk}^+ _a f_*L$
\end{itemize}
\end{lem}

\begin{proof}

Let ${\mathcal S}$ be the image of a set-theoretical section of the projection map
${\widehat A}_d\rightarrow {\widehat K}_d$ (hence of cardinality $d^{2q_f}$). Hence we have:

$$
f^{(d)}_*(L_0^{(d)})=\bigoplus_{\tau \in {\mathcal S}}\nu_d^*(f_*(L_0\otimes \tau))
$$

In characteristic 0, the pullback of a semistable vector bundle is again semistable, of the same rank and the degree is multiplied by the degree of the map. So we deduce that the  $HN$ filtration of $f^{(d)}_*(L_0^{(d)})$ has the same length $l$, slopes ${\mu_i}^{(d)}=d^{2{\overline b}}{\mu_i}$, degrees $d_i^{(d)}=d^{2q}d_i$ and ranks $r_i^{(d)}=d^{2q_f}r_i$.

\end{proof}

\bigskip
\bigskip

We proceed now to compare the continuous degree and the continuous rank of a line bundle $L$. This result will be used in last section to relate Slope Inequalities with Clifford-Severi inequalities.

\bigskip

\begin{prop}\label{gradorango} Under the usual hypothesis we have:
\begin{itemize}
\item [(i)] if $b=0$ then ${\deg^+ _a f_*L}=h^0_a(X,L)-{\rm rk}^+_a f_*L<h^0_a(X,L)$;
\item [(ii)] if ${\overline b}\geq 1$ then ${\deg^+ _a f_*L}\geq h^0_a(X,L)$.
\end{itemize}
\end{prop}

\begin{proof} (i)  If $b=0$ and $\alpha\in U_0$ we have $$h^0_a(X,L)=h^0(B,f_*(L\otimes \alpha))=h^0(B,f_*(L\otimes \alpha)^+)={\deg^+ _a f_*L}+ {\rm rk}^+_a f_*L<h^0_a(X,L)$$

\noindent by Riemann-Roch formula.

\noindent (ii) Assume now that $b\geq 1$. Let ${\mathcal G}$ be a rank $r$ semistable vector bundle on $B$ of nonnegative slope $\mu$. Recall that $h^1(B,{\mathcal G})=0$, provided ${\mu}\geq 2b-2$. If ${\mu}\geq 2b-2$ then, by Riemann-Roch formula, we have $\deg {\mathcal G}=h^0(B,{\mathcal G})+r(b-1)\geq h^0(B,{\mathcal G})$. If ${\mu}< 2b-2$ then we can apply Proposition C in \cite{H} and obtain that $\deg {\mathcal G}\geq 2h^0(B,{\mathcal G})-2r$. So, in any case, we can conclude that $\deg {\mathcal G}\geq h^0(B,{\mathcal G})-2r$.

Using this result, and arguing by induction on the length of the Harder-Narashiman filtration of any nef vector bundle ${\mathcal G}$, we obtain exactly the same result for any nef vector bundle on $B$.

Now, assume ${\overline b}\geq 1$ (which implies $b\geq 1$). Apply the above inequalities to $({\mathcal E}^{(d)})^+$. We obtain that

$$d^{2q} \deg ({\mathcal E})^+= \deg ({\mathcal E}^{(d)})^+\geq h^0(B,({\mathcal E}^{(d)})^+)-2r_m^{(d)}=d^{2q}h^0_a(X,L)-2d^{2q_f}r_m.$$

Since ${\overline b}\geq 1$, we have that $q>q_f$ and the result follows taking limit for $d\rightarrow +\infty$.

\end{proof}

\medskip

\begin{rem} If $b\geq 1$ but ${\overline b}=0$, then any inequality between the continuous degree and the continuous rank may occur. Observe that this case does not happen if $a={\alb}_X$ since then $b={\overline b}$.

\end{rem}

\section{The continuous Xiao's method}

We maintain the notations of the previous sections. In order to construct a continuous version of the classical Xiao's method, we first recall the concept of {\it continuous moving part} of a line bundle (see \cite{B}, \cite{BPS3} and \cite{J} for an alternative construction). Given a line bundle $L$ on $X$ such that $h^0_a(X,L)>0$, up to a blow-up, there is a decomposition $L=M+W$ such that, for any $d>>0$ and divisible and any general $\alpha \in {\widehat  A}$, $M^{(d)}$ is base point free and it is the moving part of $|L^{(d)}\otimes \alpha|$ and $W^{(d)}$ is its fixed divisor. Following \cite{B}, $M$ and $W$ are called the {\it continuous moving part} and {\it continuous fixed part} of $L$, respectively. According to the notation of \cite{J}, we will set $L_c=M$ for the continuous moving part.

More concretely, $L_c$ can be defined as follows. Let $V$ be any non empty open subset of $V=\{ \alpha \in {\rm Pic}^0(A)\, |\, h^0(X,L\otimes \alpha)=h^0_a(X,L)\}$. Consider the continuous evaluation map

$${\rm ev}_V: \, \bigoplus_{\alpha \in V}H^0(X,L\otimes \alpha)\otimes \alpha ^{-1} \longrightarrow L.$$

Up to a blow up $\sigma: {\overline X}\longrightarrow X$, we have that $L_c:={\rm Im} ({\overline {\rm ev}}_V)$ is a subline bundle of ${\overline L}=\sigma ^*(L)$. It turns out that $L_c$ is continuously globally generated, so it is a nef line bundle. By abuse of notation, if no confusion arises, we sometimes consider $L_c={\rm Im} ({\rm ev}_V)$ as a rank 1 subsheaf of $L$ in $X$.

\medskip

Let us consider now the natural map $f^*f_*L\rightarrow L$. It is a well known fact that its image $L'$ is a rank 1 torsion free sheaf such that, on the general fibre $F$ of $f$ it coincides with the one generated by the global sections $H^0(F,L_{|F})\otimes {\mathcal O}_F\rightarrow L_{|F}$.

\medskip

We can construct now the basic combinatorial data of Xiao's method in the continuous setting. For any $\alpha \in {\mathcal U}_0$, consider the $HN$ filtration of $f_*(L\otimes \alpha)^+$ given by $\{{\mathcal E}_i^{\alpha}\}_{i=1,..,m}$, with ranks, degrees and slopes $(r_i,d_i,\mu_i)$ not depending on $\alpha$. Fix $\alpha _0$ very general  and $L_0=L\otimes \alpha _0$ such that $\alpha_0+\tau \in {\mathcal U}_0$ for any $\tau \in {\widehat A}_{{\rm tors}}$.

For any $i\in\{i,...,m\}$ consider the continuous relative evaluation map

$$
\psi_i: \bigoplus_{\alpha \in {\mathcal U}_0}f^*({\mathcal E}_i^{\alpha})\otimes (\alpha_0-\alpha)\hookrightarrow \bigoplus_{\alpha \in {\mathcal U}_0} f^*(f_*(L\otimes \alpha))\otimes (\alpha_0-\alpha)\longrightarrow L_0.
$$

Consider a blow-up $\sigma: {\overline X} \longrightarrow X$, common for all $i$, such that the corresponding relative evaluation map ${\overline \psi_i}$ verifies that ${\rm Im} {\overline \psi} _i=T_i$ is a line bundle inside ${\overline L_0}=\sigma^*L_0$. We define $Z_i\geq 0$ to be ${\overline L_0}=T_i+Z_i$.  We define  ${\overline f}$, ${\overline F}$ and ${\overline a}$ coherently, composing with $\sigma$.


\begin{lem}\label{todo}
\begin{itemize}
\item [(i)] $L_c\leq T_m$.
\item [(ii)] For all $i=1,...,m$, ${h^0}_{{\overline a}_{|{\overline F}}}({\overline F},{T_i}_{|{\overline F}})\geq r_i$.
\item [(iii)] For all $i=1,...,m$, $N_i:=T_i-\mu_i {\overline F}$ is a nef $\mathbb{Q}$-divisor.
\item [(iv)] $f_*(T_m)^+=f_*(L_0)^+.$
\end{itemize}
\end{lem}

\begin{proof}  (i) Consider the natural composition

$$
\kappa: \bigoplus_{\alpha \in {\mathcal U}_0}H^0(X,L\otimes \alpha)\otimes (\alpha_0-\alpha)\hookrightarrow \bigoplus_{\alpha \in {\mathcal U}_0} f^*(f_*(L\otimes \alpha))\otimes (\alpha_0-\alpha)\longrightarrow L_0.
$$

Since $H^0(X,L\otimes \alpha)=H^0(B,{\mathcal E}_m^{\alpha})$, we have that $\kappa$ factorizes through ${\overline \psi}_m$ and hence
$$
L_c={\rm Im}\kappa \subseteq {\rm Im}{\overline \psi_m}=:T_m.
$$

\medskip

\noindent (ii) Fix an $\alpha \in {\mathcal U}_0$. Taking direct image by $\overline f$ in the chain of maps ${\overline f}^*{\mathcal E}_i^{\alpha}\rightarrow T_i \hookrightarrow L_0$ we obtain

$${\mathcal E}_i^{\alpha}\rightarrow {\overline f}_*T_i\otimes (\alpha -\alpha_0) \hookrightarrow {\overline f}_*({\overline L}\otimes \alpha)=f_*(L\otimes \alpha).$$

Since the composition is the natural inclusion ${\mathcal E}_i^{\alpha} \hookrightarrow f_*(L_0)$ we have that the map ${\mathcal E}_i^{\alpha}\rightarrow {\overline f}_*T_i$ is injective and so their ranks verify $r_i \leq {h^0}_{{\overline a}_{|{\overline F}}}({\overline F},{T_i}_{|{\overline F}})$.

\medskip

\noindent (iii) For $\alpha \in U_0$, $h^0({\overline F}, T_i\otimes (\alpha-\alpha_0))$ are constant, and the set $\Gamma :=\{\alpha_0 + \tau\}_{\tau \in {\widehat A}_{\rm tors}}$ is dense in ${\mathcal U}_0$. Hence we also have that

$$T_i={\rm Im}\left(\bigoplus_{\alpha \in {\Gamma}}{\overline f}^*({\mathcal E}_i^{\alpha})\otimes (\alpha_0-\alpha)\longrightarrow {\overline L}_0\right).$$

Hence, by the Noetherian property, $T_i$ must be the image of a {\it finite} sum of such maps associated to torsion $\tau_1,...,\tau_e$, of orders $d_1,...,d_e$. Let $d=d_1...d_e$ and consider $\nu_d$ the multiplication map by $d$. Then

$$
\nu_d^*(T_i)=\pi_i^*({\mathcal O}_{\mathbb{P}^{(d)}_i}(1)),
$$

\noindent where
$$
\pi_i:{\overline X}^{(d)}\longrightarrow \mathbb{P}^{(d)}_i:=\mathbb{P}_{B^{(d)}}(\oplus_{\tau \in {\widehat A}_d}\nu_d^*{\mathcal E}_i^{\alpha_0+\tau}).
$$

If $R$ is a general fibre of the projective bundle, then ${\mathcal O}_{\mathbb{P}^{(d)}_i}(1)-\mu_iR$ is nef by \cite{M}, then $\nu_d^*(T_i)-\mu_iF^{(d)}=(T_i-\mu_i{\overline F})^{(d)}$ is nef and so is $T_i-\mu_i{\overline F}$. This fact is a continuous version of (\cite{X} Lemma 4).

\medskip

\noindent (iv) By construction we have $f_*(L_0)^+\subseteq f_*(T_m)\subseteq f_*(L_0)$. Then take the $+$-part, which is the biggest nef subbundle inside them (see Lemma \ref{remarkpositivepart}).
\end{proof}

\bigskip

By construction, we have that $T_i\otimes \beta + Z_i ={\overline L_0}\otimes \beta$, for any $\beta \in -\alpha_0+{\mathcal U}_0$.

Denote by $P_i:={T_i}_{|{\overline F}}$. The class $P_i\otimes {{\overline a}_{|{\overline F}}^*\beta}$ is then parametrized by a non empty open set of ${\widehat K}$. Hence, after modifying $\alpha_0$, if necessary, by (ii) of the previous lemma, we may assume that

$$
h^0({\overline F},P_i)=h^0_{{\overline a}_{|{\overline F}}}({\overline F},P_i)\geq r_i.
$$

\medskip

We have constructed on ${\overline X}$ sequences

$$N_1 \leq N_2 \leq ...\leq N_m  \leq {\overline L_0}$$
$$Z_1 \geq Z_2 \geq ... \geq Z_m \geq  0$$
$$\mu_1 \geq \mu_2 \geq ... \geq \mu_m \geq 0$$

\noindent of nef $\mathbb{Q}-$divisors $N_i$, effective divisors $Z_i$ and rational numbers $\mu_i$ such that, for any $i$, ${\overline L_0}=N_i+\mu_i{\overline F}+Z_i$.

We can define coherently

$$(N_{m+1},Z_{m+1},\mu_{m+1})=(T_m,Z_m,0)).$$



Under these assumptions we can apply Xiao's Lemma (see \cite{K} and \cite{O} Lemma 1.2), to produce a broad set of inequalities. Fix $r\leq n$ and consider an ordered partition of the set $\{ \, 1,...,m \, \}$ given by subsets $I_s$, for $s=1,...,r-1$ and such that $I_{r-1}\neq \emptyset$ (some of the other $I_i$ may be empty).

Define decreasingly, for $s=1,...,r-1$

$$b_s=\{
           \begin{array}{cc}
             {\rm min} I_s & {\rm if} \,\, I_s \neq \emptyset \\
             b_{s+1} & {\rm otherwise} \\
           \end{array}$$

\bigskip

 Then we have that, for any $Q_1,...,Q_{n-r}$ nef $\mathbb{Q}$-Cartier divisors the following inequality holds:

\begin{equation}\label{Xiaogeneral}
Q_1...Q_{n-r}\left[N^r_{m+1}-(\sum _{s=r-1}^{1} (\prod_{k>s}P_{b_k})\sum_{i \in I_s}(\sum_{l=0}^{s}P_i^{s-l}P_{i+1}^{l})(\mu_i-\mu_{i+1}))\right]\geq 0.
\end{equation}

\bigskip

\noindent For further reference, if we take $r=n$ and $I_{n-1}=\{1,...,m\}$, we obtain

\begin{equation}\label{xiaofinito}
{\rm vol}(L)={\rm vol} (L_0)={\rm vol}({\overline L_0})\geq N_{m+1}^n\geq \sum_{i^1}^m \left[P_{i+1}^{n-1}+P_{i+1}^{n-2}P_i+...+P_i^{n-1}\right](\mu_i-\mu_{i+1}).
\end{equation}

The use of this inequality is based on the following. Observe that

$$
{\rm deg}_a^+f_*(L)={\rm deg}({\mathcal E}_m^{\alpha_0})=\sum_{i=1}^{m}r_i(\mu_i-\mu_{i+1}).
$$

Hence, to obtain a lower bound of ${\rm vol}(L)$ in terms of ${\rm deg}_a^*f_*(L)$ it is enough to find lower bounds of

$$
\delta_i=P_{i+1}^{n-1}+P_{i+1}^{n-2}P_i+...+P_i^{n-1}
$$

\noindent in terms of $r_i$. Recall that, by construction, $P_i$ is a linear system on ${\overline F}$ such that $
h^0({\overline F},P_i)=h^0_{{\overline a}_{|{\overline F}}}({\overline F},P_i)\geq r_i.$ For later reference observe that, since $P_i\leq P_{i+1}$ and they are nef, we have

\begin{equation}\label{delta}
\delta_i\geq P_{i+1}^{n-1}+(n-1)P_i^{n-1}.
\end{equation}


\bigskip
\bigskip

\section{Slope inequalities and equalities}


\bigskip

We will consider triplets $(X,L,a)$ where $X$ is an irregular variety, $L$ a big line bundle on $X$ and $a:X\longrightarrow A$ a nontrivial generating map to an abelian variety. We will always assume that $h^0_a(X,L)\neq 0$. We can also consider {\it fibred} triplets, $(X,L,a,f)$ where $f:X:\longrightarrow B$ is an {\it irregular} fibration onto a smooth curve of genus $b$, as in Section 2.

\begin{defn}
\begin{itemize}
\item [(i)] We define the {\it Clifford-Severi slope} of $(X,L,a)$ as

$$
\lambda(L,a)=\frac{{\rm vol} L}{h^0_a(X,L)}\in(0,+\infty).
$$

\item [(ii)] We define the {\it slope of the fibred triplet} as

$$
s(f,L,a)=\frac{{\rm vol} L}{{\rm deg}_a^+(L)}\in (0,+\infty].
$$

\end{itemize}
\end{defn}

\begin{rem} We are interested in relating both slopes. For a fixed irregularly fibred triplet, in some cases we have a clear inequality between them, given by Proposition \ref{gradorango}. For example, $s(f,L,a)\geq \lambda (L,a)$ when $b=0$, but this does not hold for ${\overline b}\geq 1$.

As we will see, both slopes are closely related when we consider their minimum value among triplets on a good family of triplets ${\mathcal F}$. {\it Clifford-Severi inequalities} and {\it Slope inequalities} arise when we study such minima, respectively.
\end{rem}

\smallskip

We can impose several conditions to define families ${\mathcal F}$ for which we know non-trivial Clifford-Severi inequalities in the literature.

\begin{defn}\label{conditions}

\begin{itemize}
\item[(i)] \underline{Condition $P_r$ ($0<r\leq +\infty$)}: $L$ is a numerically $r$-subcanonical line bundle (see Definition 5.1 in \cite{BPS3}). If $r=1$ we will say that $L$ is {\it numerically subcanonical}. If $r=+ \infty$, then $L$ is any line bundle.
\item[(ii)] \underline{Condition $Q_c$}: $\dim X- \dim a(X)\leq c$. The case $c=0$ means that $X$ is of maximal $a$-dimension.
\item[(iii)] \underline{Conditions on the degree of $a$} when $X$ is of maximal $a$-dimension. $D_1$: $a$ is of degree 1. $D_{\neq 2}$: $a$ is not composed with an involution.
\item[(iv)] \underline{Condition $G$}: $X$ is of general type.
\end{itemize}
\end{defn}


\medskip

\begin{defn}\label{goodclass} Let  ${\mathcal F}$ be a class of triplets $(X,L,a)$ such that $a:X\longrightarrow A$ is a generating map to an abelian variety and $h^0_a(X,L)>0$. We say that the class is {\it good} if it verifies the following conditions:
\begin{itemize}
\item[(1)] If $(X,L,a)\in {\mathcal F}$, then $({\overline X},\sigma ^*L,a\circ \sigma)\in {\mathcal F}$, where $\sigma:{\overline X}\longrightarrow X$ is a birational morphism.
\item[(2)] If $(X,L,a)\in {\mathcal F}$, then $(X^{(d)},L^{(d)},a_d)\in {\mathcal F}$.
\item[(3)] If $(X,L,a)\in {\mathcal F}$, then $(X,L',a)\in {\mathcal F}$ for $L'\leq L$ such that $h^0_a(X,L')>0$.
\item[(4)] If $(X,L,a)\in {\mathcal F}$, then $(M,L_{|M},a_{|M})\in {\mathcal F}$, for a general smooth $M$ in any base point free linear system in $X$.

\bigskip

\end{itemize}
\end{defn}

\bigskip

\begin{rem} Any combination of conditions $P_r$, $Q_c$, $D$ and $G$ define good classes.
\end{rem}
\bigskip

\begin{defn} Given a good class ${\mathcal F}$ and a triplet $(X,L,a)$ (not necessarily in ${\mathcal F}$), with an irregular fibration $f: X \longrightarrow B$, we say that {\it the fibration $f$ if of type} ${\mathcal F}$ if $(F,L_{|F},a_{|F})$ is in ${\mathcal F}$, for a general fibre $F$ of $f$.
\end{defn}

\begin{rem} Observe that this is equivalent to say that $(X,L,a)$ is in ${\mathcal F}$ in most cases, but not always. For example, when $B=\mathbb{P}^1$, $X$ being of maximal Albanese dimension implies $f$ is of maximal Albanese dimension, but the viceversa is not true.
\end{rem}

\bigskip

\begin{defn} Given a class ${\mathcal F}$ of triplets we define:

\begin{itemize}
\item[(i)] The minimal {\it Clifford-Severi} slopes for ${\mathcal F}$:

$${\lambda}_{\mathcal F}(n)={\rm inf}\{\lambda (L,a)|\,(X,L,a)\in {\mathcal F},\,n=\dim X\}.$$

\item[(ii)] The minimal {\it irregular Fibration} slopes for ${\mathcal F}$:

$$s_{\mathcal F}(n)={\rm inf}\{s (L,f,a)\,|\, f \,{\rm irregular} \,\,{\rm of}\,\,{\rm type}\,\, {\mathcal F},\,n=\dim X\}.$$

\item[(iii)] Expressions of type ${\lambda}_{\mathcal F}(n)\geq \lambda_n$ or $s_{\mathcal F}(n)\geq \lambda_n$ will be called  {\it Clifford-Severi inequalities} and {\it Slope inequalities}, respectively, for the class $\mathcal F$.

\end{itemize}


\end{defn}
\bigskip

\begin{rem}\label{CS} For varieties of maximal $a$-dimension, a list of Clifford-Severi Inequalities (and Equalities) we can find in the literature for some good classes is (cf. \cite{B}, \cite{BPS3}):

\begin{itemize}
\item If ${\mathcal F}$ is defined by $P_{\infty}$ then ${\lambda}_{\mathcal F}(n)= n!$.
\item If ${\mathcal F}$ is defined by $P_1$ then ${\lambda}_{\mathcal F}(n)= 2 \,n!$.
\item If ${\mathcal F}$ is defined by $P_r$ and $D_{\neq 2}$, then ${\lambda}_{\mathcal F}(n)\geq \delta (r)\, n!$.
\item If ${\mathcal F}$ is defined by $P_r$ and $D_1$, then ${\lambda}_{\mathcal F}(n)\geq {\overline \delta (r)}\,n!$.
\item If ${\mathcal F}$ is defined by $P_1$ and $D_{\neq 2}$, then ${\lambda}_{\mathcal F}(n)\geq \frac{9}{4}\,n!$.
\item If ${\mathcal F}$ is deined by $P_1$, $G$ and $D_1$, then ${\lambda}_{\mathcal F}(n)\geq \frac{5}{2}\, n!$.
\end{itemize}

\noindent where $\delta(r)=2$ if $r\leq 1$ and $2r/(2r-1)$ otherwise and ${\overline \delta} (r)$ has a piecewise similar expression (see Section 6 in \cite{BPS3}).
\end{rem}

\bigskip

\begin{rem}\label{csequality}
In \cite{BPS4}, the authors characterize varieties attaining the equality in the first two cases, assuming $h^0_a(X,L)\neq 0$:

\begin{itemize}
\item If ${\rm vol} (L)=n!\, h^0_a(X,L)$, then $\dim A=\dim X$ and the map $a$ is birrationally étale.
\item If, $n\geq 2$, $L$ is numerically subcanonical and ${\rm vol} (L)=2\,n!\,h^0_a(X,L)$, then $\dim A=\dim X$ and the map $a$ is birrationally a double cover composed with an étale map on $A$.
\item If $n=1$, $L$ is numerically subcanonical and ${\rm vol} (L)=2\,h^0_a(X,L)$, then, either $L=K_X$ or $a$ is a double cover of an elliptic curve $A=E$, composed with an étale map, and $L\in a^*{\rm Pic} (E)$.
\end{itemize}
\end{rem}

\bigskip

Given a good class $\mathcal F$, we can consider similar problems for irregular fibrations of type $\mathcal F$, and would like to find lower bounds for $s_{\mathcal F}$ (see \cite{BS} for a survey on this topic). Although there are not general slope inequalities in higher dimension, recently Hu and Zhang (\cite{HZ}) obtained a set of such slope inequalities for the relative canonical bundle $K_f$ by a involved direct computation. Here, we are going to give a broad generalization and see that for any good class ${\mathcal F}$, Clifford-Severi inequalities and Slope inequalities are equivalent, in some strong sense. The equivalence is not dimension to dimension, but with a zigzag argument, and allows to obtain Slope inequalities automatically. The proof is completely straightforward using the continuous Xiao's method constructed in Section 3, and allows to characterize easily the equality cases. More concretely:

\smallskip

\begin{thm}\label{metateorema}
Let ${\mathcal F}$ be a good class of triplets of maximal $a$-dimension. Then
\begin{itemize}
\item [(i)]  $s_{\mathcal F}(n)\geq \,n\, {\lambda}_{\mathcal F}(n-1)$.

\smallskip

\item [(ii)] ${\lambda}_{\mathcal F}(n)\geq s_{\mathcal F}(n)$.

\smallskip

\item [(iii)] If ${\lambda}_{\mathcal F}(n-1)\geq \lambda$, and for a fibred triplet $(X,L,a)\in {\mathcal F}$ we have $s(f,L,a)= n\,\lambda$, then
        \begin{itemize}
        \item [(1)] For a general $\alpha \in {\widehat A}$, the vector bundle $f_*(L\otimes \alpha)^+$ is semistable.
        \item [(2)] There exists $L'_F\leq L_{|F}$ such that $\lambda={\lambda}_{\mathcal F}(n-1)=\lambda ({L'_F},a_{|F})$.
        \end{itemize}
        
\smallskip
        
\item [(iv)] If in (iii) we have that $L$ is nef, we can also conclude that
        \begin{itemize}
        \item [(3)] For a general $\alpha \in {\widehat A}$, the vector bundle $f_*(L\otimes \alpha)$ is nef (and hence semistable).
        \item [(4)] $\lambda={\lambda}_{\mathcal F}(n-1)=\lambda (L_{|F},a_{|F})$.

        \end{itemize}

\end{itemize}
\end{thm}

\bigskip

\begin{rem} As pointed out in the introduction (see Corollary A) this theorem allows to obtain, in an automatic way, higher dimensional Slope and Clifford-Severi inequalities for good classes of varieties of maximal $a$-dimension, just knowing a Clifford-Severi inequality with  $\lambda_{m_0}$ for a certain low dimension $m_0$ (typically $m_0=1$ or $2$). In particular, this provides a new proof of all known higher dimensional Clifford-Severi inequalities given in Remark \ref{CS}. This strategy is also used in a different way in the proofs by induction on the dimension given in \cite{BPS2}, just integrating the volume and continuous rank functions. In some sense, we can consider that the technique of the Continuous Xiao's method is a discretization of the integral of such functions used there.
\end{rem}

\bigskip

\begin{proof}

(i) Apply continuous Xiao's method to the fibration $f$. Let $\lambda:=\lambda_{\mathcal F}(n-1)$. Since $(F,L_{|F},a_{|F})$ are in ${\mathcal F}$, so are
$({\overline F},P_i,{\overline a}_{|{\overline F}})$, by conditions (1),(2) and (3) in Definition \ref{goodclass}. So we obtain that
        $$
        P_i^{n-1}\geq \lambda \, h^0_{{\overline a}_{|{\overline F}}}({\overline F},P_i)\geq \lambda \, r_i.
        $$

Then, by (\ref{delta}) we have

\begin{equation}\label{casolimite}
\delta_i=\left[P_{i+1}^{n-1}+P_{i+1}^{n-2}P_i+...+P_i^{n-1}\right]\geq P_{i+1}^{n-1}+(n-1)P_i^{n-1}\geq \lambda (r_{i+1}+(n-1)r_i)\geq n\lambda r_i.
\end{equation}

\noindent On the other hand we have ${\deg}^+_af_*L={\deg} (f_*L_0)^+=\sum_{i=1}^mr_i(\mu_i-\mu_{i+1})$ and so by applying (\ref{xiaofinito}) we obtain

\begin{equation}\label{volumenesiguales}
{\rm vol} (L)\geq N_{m+1}^n\geq n\, \lambda \,{\deg^+_a f_*(L)}.
\end{equation}

\bigskip

\noindent (ii) This results follows directly by applying Pardini's trick (\cite{P}). We sketch the proof for benefit of the reader. Take a triple $(X,L,a) \in {\mathcal F}$ of dimension $n$ and let $\lambda=s_{\mathcal F}(n)$. We need to prove that $\lambda (L,a)\geq \lambda$.

For this, take a fixed very ample line bundle $H$ on $A$ and a multiplication map $\nu_d$. Consider the following diagram:

\begin{equation}\label{Pardinistrick}
\xymatrix{
Y_d\ar[d]_{f_d} \ar[rd]_{{\widetilde a}_d}\ar[r]_{\sigma}&X^{(d)}\ar[d]_{a_d}\ar[r]_{\nu_d} &X\ar[d]^a\\
\mathbb{P}^1&A\ar[r]_{\nu_d}&A}
\end{equation}

\noindent constructed as follows. Take two general smooth members of $|a_d^*(H)|$  and let $\sigma$ be the blow-up of the base locus of the pencil generated by them. $f_d$ is the resulting fibration and call $M_d$ a general fibre (observe that $M_d\in |a_d^*(H)|$). Let $M\in |a^*(H)|$.

Let $L^{(d)}=\nu_d^*(L)$ and $L_d=\sigma^*L^{(d)}$. We have the following numerical relations (see \cite{BPS3} for details):

$${\rm vol}(L_d)={\rm vol}(L^{(d)})=d^{2q}{\rm vol}(L),$$
$$M_d\equiv \frac{1}{d^2}M^{(d)},$$
$$h^0_{a_d}(X^{(d)},L^{(d)})=d^{2q}h^0_a(X,L).$$

By Proposition \ref{gradorango} (i), we have that

$${\rm deg}_{{\widetilde a}_d}^+(f_d)_*(L_d)=h^0_{{\widetilde a}_d}(Y_d,L_d)-{\rm rk}_{{\widetilde a}_d}^+(f_d)_*(L_d).$$

Observe now that

$$h^0_{{\widetilde a}_d}(Y_d,L_d)=h^0_{a_d}(X^{(d)},L^{(d)})=d^{2q}h^0_a(X,L),$$
$${\rm rk}_{{\widetilde a}_d}^+(f_d)_*(L_d)\leq {\rm rk}_{{\widetilde a}_d}(f_d)_*(L_d)=h^0_{a_d}(M_d,L^{(d)})$$

In order to bound the last term we apply the general Clifford-Severi inequality for irregular varieties given in \cite{B} (Main Theorem): $h^0_{a_d}(M_d,L^{(d)})\leq {\rm vol}_{M_d}(L^{(d)})$. Finally, observe that ${\rm vol}_{M_d}(L^{(d)})\leq \frac{1}{d^{2q-2}}{\rm vol}_M(L)$. Indeed, both volumes agree if $L$ is nef since then the volume is computed by an intersection product and $M_d\equiv \frac{1}{d^2}M^{(d)}$. When $L$ is not nef, we approximate the computation of both volumes by nef $Q$-line bundles using Fujita Approximation theorem.

Now, taking $d\rightarrow +\infty$ we obtain that

$$\lambda (L,a)\geq s(f_d,L_d,{\widetilde a}_d).$$

Finally, observe that $s(f_d,L_d,{\widetilde a}_d)\geq \lambda$, since $f_d$ is a fibration of type ${\mathcal F}$. Indeed, observe that, by Definition \ref{goodclass}, if $(X,L,a)$ is in ${\mathcal F}$, so is $(X^{(d)},L^{(d)},a_d)$ and hence $(M_d,L^{(d)}_{|M_d},({a_d})_{|M_d})$ and $(M_d,({L_d})_{|M_d},({\widetilde a}_d)_{|M_d})$ are also in ${\mathcal F}$. Hence $f_d$ is of type ${\mathcal F}$, and it is an irregular fibration since $b=0$.

\bigskip

\noindent (iii) In case of equality in (i), for all $i$ we have $P_i^{n-1}= \lambda\,h^0_{{\overline a}_{|{\overline F}}}({\overline F},P_i)=\lambda\,r_i$, and hence, for all $i$ the equality

$$h^0_{{\overline a}_{|{\overline F}}}({\overline F},P_i)=r_i$$

\noindent holds.

Moreover, by (\ref{casolimite}) we have that $r_{i+1}=r_i$ holds for all $i$. Then, necessarily $m=1$ and so for a general $\alpha \in {\widehat A}$ the vector bundle $f_*(L\otimes \alpha)^+$ is semistable.

In particular,  ${T_1}_{|{\overline F}}^{n-1}=P_1^{n-1}=\lambda \,h^0_{{\overline a}_{|{\overline F}}}({\overline F},P_1))$ and we can take $L'_F:=T_1$.


\bigskip

\noindent (iv) Another consequence of the equality is that ${\vol L}=N_{m+1}^n$ (see \ref{volumenesiguales}). Since $L$ is nef by hyothesis and big, we can apply \cite{FKL}, Theorem A, and obtain that $Z_1={\overline L}_0-N_{m+1}=0$. Then $N_{m+1}=T_m=T_1$ implies that

$${\rm rank}f_*(L_0)^+=r_1=h^0_{{\overline a}_{|{\overline F}}}({\overline F},P_1)=h^0_{{\overline a}_{|{\overline F}}}({\overline F},{L_0}_{|{\overline F}})=h^0_{a_{|F}}(F,L_0)={\rm rank}f_*(L_0).$$

Hence we conclude that $f_*(L_0)=f_*(L_0)^+$ and so it is nef (which is (3)), and that $(L_{\F},a_{|F})$ verifies the Clifford-Severi equality (which is (4)).

\end{proof}

\bigskip

\begin{rem} Theorem \ref{metateorema} holds under the condition of maximal $a$-dimension. However, some of the statements hold for classes of non maximal $a$-dimension. For example, observe that for inequality (i) we only need properties (1),(2) and (3) of a good class. Nevertheless, in this form, the statement is trivial since general Clifford-Severi inequalities for such varieties have $\lambda_{\mathcal F}(r)=0$ for all $r$. The reason for this is that we need to impose some kind of bigness on line bundles in order to have non trivial inequalities. But bigness do not verify condicion (3) in the definition of a good class, necessary to use Xiao's method.

The fibration slope inequalities for non maximal $a$-dimension varieties are not equivalent to Clifford-Severi ones, and need a slightly different strategy. This will be considered in a forthcoming paper (\cite{B2}.
\end{rem}

\bigskip

From the previous meta-theorem, applied to the cases listed in Remark \ref{CS} we automatically obtain a series of slope inequalities, and the corresponding characterization of equalities.
\bigskip

\begin{cor} Let $(X,L,a)$ a triple with a fibration $f: X \longrightarrow B$ onto a smooth curve. Assume that $h^0_a(X,L)\neq 0$ and that $f$ is of maximal $a$-dimension.
\begin{itemize}
\item If $L_{|F}$ is numerically $r$-subcanonical and $a_{|F}$ is not composed with an involution, then ${\rm vol}(L)\geq \, \delta (r) n!\, {\rm deg}_a^+f_*L.$
\item If $L_{|F}$ is numerically $r$-subcanonical and $a_{|F}$ is of degree 1, then ${\rm vol}(L)\geq \, {\overline \delta}(r)n!\, {\rm deg}_a^+f_*L.$
\item If $L_{|F}$ is numerically subcanonical and $a_{|F}$ is not composed with an involution, then ${\rm vol}(L)\geq \, \frac{9}{4}n!\, {\rm deg}_a^+f_*L.$
\item If $F$ is of general type, $L_{|F}$ is numerically subcanonical and $a_{|F}$ is of degree 1, then ${\rm vol}(L)\geq \, \frac{5}{2}n!\, {\rm deg}_a^+f_*L.$
\end{itemize}

In all the previous cases, if equality occurs and $L$ is nef, then $f_*(L\otimes \alpha)$ is nef and semistable, for general $\alpha$, and $(F,L_{|F},a_{|F})$ verifies the corresponding Clifford-Severi equality.
\end{cor}

\bigskip

There are two cases of special relevance where we can say even more:

\bigskip

\begin{cor}\label{cor1}
  Let $X$ be a variety of dimension $n$ and $a:X\longrightarrow A$ be a generating map to an abelian variety of dimension $q$. Let $f:X\longrightarrow B$ a fibration of maximal $a-$dimension and let $b=g(B)$. Let $L$ be line bundle on $X$ such that $h^0_a(X,L)\neq 0$. Then:
  \begin{itemize}
  \item [(i)] ${\rm vol} (L)\geq n!\,{\deg^+_a f_*(L)}.$
  \item [(ii)] If equality holds, then:
      \begin{itemize}
       \item [$(1)$] For a general $\alpha \in {\widehat A}$, the vector bundle $f_*(L\otimes \alpha)^+$ is semistable. If $L$ is nef then $f_*(L\otimes \alpha)$ is nef and semistable.
      \item [$(2)$] There exist $L'_F \leq L_{|F}$ on the Severi line: ${\rm vol}(L'_F)=(n-1)!h^0_{a_{|F}}(F,L'_F).$ If $L$ is nef, we can take $L'_F=L_{|F}$.
      \item [$(3)$] Assume that $h^0_{a_{|F}}(F,L_{|F})\neq 0$. Then we have $q-{\overline b}=n-1$ and the map $a_{|F}:F \longrightarrow K$ is, birationaly, an étale map on $K$.
      \end{itemize}
  \end{itemize}
\end{cor}

\begin{proof} Just consider Theorem \ref{metateorema}, and Remarks \ref{CS} and \ref{csequality}.
\end{proof}

\bigskip

\begin{cor}\label{cor2}
  Let $X$ be a variety of dimension $n$ and $a:X\longrightarrow A$ be a generating map to an abelian variety of dimension $q$. Let $f:X\longrightarrow B$ a fibration of maximal $a-$dimension and let $b=g(B)$. Let $L$ be a numerically subcanonical line bundle on $X$ such that $h^0_a(X,L)\neq 0$. Then:
  \begin{itemize}
  \item [(i)] ${\rm vol} (L)\geq 2\,n!\,{\deg^+_a f_*(L)}.$
  \item [(ii)] If equality holds, then:
      \begin{itemize}
       \item [$(1)$] For a general $\alpha \in {\widehat A}$, the vector bundle $f_*(L\otimes \alpha)^+$ is semistable. If $L$ is nef then $f_*(L\otimes \alpha)$ is nef and semistable.
      \item [$(2)$] There exist $L'_F \leq L_{\F}$ on the Severi line: ${\rm vol}(L'_F)=2(n-1)!h^0_{a_{|F}}(F,L'_F).$ If $L$ is nef, we can take $L'_F=L_{|F}$.


      \item [$(3)$] Assume that $h^0_{a_{|F}}(F,L_{|F})\neq 0$. Except in the case $n=2$ and $L_{|F}=K_F$, we have $q-{\overline b}=n-1$ and the map $a_{|F}:F \longrightarrow K$ is, birationaly, a double cover composed with an étale map on $K$.
     \end{itemize}
  \end{itemize}
\end{cor}

\begin{proof}

We only have to prove $(ii)-(3)$. By Theorem \ref{metateorema} (iii), we have that $P_1\leq L_{|F}$ is a numerically subcanonical divisor verifying $P_1^{n-1}= 2\,(n-1)!\,h^0_{a_{|F}}(F, P_1)$. By Remark \ref{csequality}, if $n\geq 3$, $q-{\overline b}={\rm dim}K=n-1\geq 2$ and so the map $a_{|F}$ is, birationally, a double cover composed with an étale map on $K$. If $n=2$ the same conclusion holds except in the case $P_1=K_F$. Since $K_F=P_1\leq L_{|F}\leq {K_X}_{|F}=K_F$, then we would have $L_{|F}=K_F$.

\end{proof}

In the particular case  $L=K_f$ (the relative canonical divisor of the fibration $f$)  and $a={\rm alb}_X$, we can obtain even more information (compare with \cite{HZ}, Theorem 1.3). This result is the higher dimension analog of Xiao's theorem for irregularly fibred surfaces of minimal slope (see Theorem 3 in \cite{X}).

\bigskip

\begin{cor}\label{corolario}
Let $X$ be a complex, projective, irregular, $\mathbb{Q}$-Gorenstein variety of dimension $n$ with at most terminal singularities. Let $f:X \longrightarrow B$ a relatively minimal fibration of maximal Albanese dimension and let $F$ be a general fibre. Then:
\begin{itemize}
\item[(i)] $K_f^n\geq 2\,n! \chi_f.$
\item[(ii)] If equality holds, then
    \begin{itemize}
    \item [$(1)$] For general $\alpha \in {\rm Pic}^0X$, the vector bundle $f_*(\omega_f\otimes \alpha)$ is nef and semistable.
    \item [$(2)$] $K_F^{n-1}=2\,(n-1)! \chi(\omega_F).$
    \item [$(3)$] Assume $\chi (\omega_F)>0$. If $n\geq 3$, then $q(F)=n-1$ and $q(X)=b+q(F)$ and the map $(alb_X)_{|F}:F\longrightarrow K$ is, birationally, a double cover composed with an étale map on $K$. Then apply Corollary \ref{cor2}.
     \end{itemize}

\end{itemize}
\end{cor}

\begin{proof}
Consider a desingularization $\sigma: {\overline X} \longrightarrow X$, ${\overline f}=\sigma \circ f$ the induced fibration adn ${\overline F}$ its general fibre. Then we have $\vol (K_{\overline X})=\vol (K_X)=K_X^n$, $\chi ({\overline X}, K_{\overline X})=\chi (X,K_X)$ and the same holds for ${\overline F}$ and $F$.
\end{proof}


\bigskip



     \begin{thebibliography}{ABCD}



\bibitem{BS} M.A.~Barja, L.~Stoppino, {\em Stability conditions and positivity of invariants of fibrations}.  Algebraic and complex geometry, 1–40, Springer Proc. Math. Stat., 71, Springer, Cham, 2014.

\bibitem{B} M.A.~Barja, {\em Generalized Clifford Severi inequality and the volume of irregular varieties}. Duke Math. J. {\bf164} (2015), no. 3, 541--568.

\bibitem{B2} M.A.~Barja, {\em Slope Inequalities for fibrations of non-maximal Albanese dimension}. In preparation.

\bibitem{BPS1} M.A.~Barja, R.~Pardini, L.~Stoppino, {\em Surfaces on the Severi line}, Journal de Math\'ematiques Pures et Appliqu\'ees, (2016), no. 5, 734--743.

\bibitem{BPS2} M.A.~Barja, R.~Pardini, L.~Stoppino, {\em The eventual paracanonical map of a variety of maximal Albanese dimension}.  Algebr. Geom. 6 (2019), no. 3, 302–311.

\bibitem{BPS3} M.A.~Barja, R.~Pardini, L.~Stoppino, {\em Linear systems on irregular varieties}. Journal of the Institute of mathematics of Jussieu, 1 (2019), p. 1-39.

\bibitem{BPS4} M.A.~Barja, R.~Pardini, L.~Stoppino, {\em Higher Dimensional Clifford-Severi equalities}. Communications in Contemporary Mathematics, (2019) doi.org/10.1142/S0219199719500792



\bibitem {Beau} A.~Beauville,
{\em Annulation du $H^1$ pour les fibr\'es en droites plats}, Complex algebraic varieties (Bayreuth, 1990), 1--15,
Lecture Notes in Math., 1507, Springer, Berlin, 1992.

















\bibitem{FKL} M.~Fulger, J.~Koll\'ar, B.~Lehmann, {\em Volume and Hilbert Functions of $\R$-divisors}, Michigan Math. J. {\bf 65} (2016), no. 2, 371--387.








\bibitem{H} G. ~Hein, D.~Ploog {\em Stability of Picard sheaves for vector bundles ov curves}. Journal of Geometry and Physics 122  (2017), 59-68.

\bibitem{HP} C.~Hacon, R.~Pardini. {\em Birational characterization of products of curves of genus 2}. Math. Res. Lett. 12 (2005), no. 1, 129–140.

\bibitem{HZ} Y. ~Hu, T.~Zhang. {\em Relative Severi inequality for fibrations of maximal Albanese dimension over curves }, arXiv: 1905.08404 [math.AG]

\bibitem{J} Z. Jiang, {\em On Severi type inequalities}, Math. Ann. (2020). doi.org/10.1007/s00208-020-02082-6

\bibitem{K} K. ~Konno {\em A  lower  bound  of  the  slope  of  trigonal  fibrations}, Internat. J. Math.7(1996), no. 1, 19–27.

\bibitem{LZ} X.~Lu, K.~Zuo {\em On the Severi type Inequalities for Irregular Surfaces}. International Mathematics Research Notices, 1 (2019), 231–248.

\bibitem{M} Y.~Miyaoka {|The Chern classes and Kodaira dimension of a minimal variety}, Algebraic Geometry, Sencai (1985). Adv. Study in Pure Math. 10 (1987), 449-476.

\bibitem{MLP} M. ~Mendes Lopes, R. ~Pardini {\em Severi type inequalities for irregular surfaces withample canonical class}, Comm. Math. Helv. 86 (2011), no. 2, 401-414.

\bibitem{O} K.~Ohno {\em Some inequalities for minimal fibrations of surfaces of general type over curves}. J. Math. Soc. Jampan 44 (1992). 75-89,





\bibitem{P} R.~Pardini, {\em The Severi inequality $K^2\geq 4\chi$ for surfaces of maximal Albanese dimension}, Invent. Math {\bf 159} {\bf 3} (2005), 669--672.

\bibitem {Sh} S.~Shatz,  {\em The decomposition and specialization of algebraic families of vector bundles}. Compositio Mathematica, 35 (1977), no. 2, 163-187.

\bibitem {S} L.~Stoppino, {\em The continuous rank function for varieties of Maximal Albanese dimension and its applications}, Geometry of Varieties and Applications (2019)


\bibitem{X}  G.~Xiao, {\em Fibered  algebraic  surfaces  with low slope}, Math. Ann.276(1987), no. 3,449–466.


\bibitem{Z} T.~Zhang, {\em Severi inequality for varieties of maximal Albanese dimension}, Math. Annalen {\bf 359} (2014), 3, 1097--1114.

\bibitem{Z2} T.~Zhang, {\em Relative Clifford inequality for varieties fibred by curves}, Journal of Differential Geometry, to appear. (  arXiv:1706.06523).


     \end{thebibliography}
     \end{document}